\documentclass{amsart}
\newcommand{\bbc}{\mathbb{C}}
\newcommand{\bbf}{\mathbb{F}}
\newcommand{\bbq}{\mathbb{Q}}\newcommand{\bbp}{\mathbb{P}}
\newcommand{\bbz}{\mathbb{Z}}
\newcommand{\bbn}{\mathbb{N}}
\newcommand{\caa}{\mathcal{A}}
\newcommand{\cad}{\mathcal{D}}
\newcommand{\caf}{\mathcal{F}}\newcommand{\cag}{\mathcal{G}}

\newcommand{\cao}{\mathcal{O}}
\newcommand{\cas}{\mathcal{S}}\newcommand{\cax}{\mathcal{X}}
\newcommand{\cay}{\mathcal{Y}}\newcommand{\caz}{\mathcal{Z}}
\newcommand{\frA}{\mathfrak{A}}\newcommand{\frB}{\mathfrak{B}}
\newcommand{\frp}{\mathfrak{p}}

\DeclareMathOperator{\Frob}{Frob}

\DeclareMathOperator{\cl}{cl} \DeclareMathOperator{\divi}{div}
\DeclareMathOperator{\Div}{Div}
\DeclareMathOperator{\et}{\text{\'et}} \DeclareMathOperator{\gal}{Gal}

\DeclareMathOperator{\ns}{NS}
\DeclareMathOperator{\ch}{CH}
\DeclareMathOperator*{\ord}{ord}
\newcommand{\ov}{\overline}
\DeclareMathOperator{\pic}{Pic} \DeclareMathOperator*{\res}{{R\acute{e}s}}
\DeclareMathOperator{\rg}{rang}
 
\DeclareMathOperator{\tr}{Tr}
\newcommand{\Dem}{D\'emonstration}

\newcommand{\un}{\underline}
\newcommand\jour{\number\day\space \space\ifcase\month\or janvier\or
f\'evrier\or mars\or avril\or mai\or juin\or
juillet\or ao\^ut\or septembre\or octubre\or novembre\or
d\'ecembre\fi\space \space\number\year}
\newtheorem{theorem}{Th\'eor\`eme}
\numberwithin{theorem}{section} \theoremstyle{plain}

\newtheorem*{conjTa}{Conjecture T}
\newtheorem*{conjMan}{Conjecture M$_{\text{an}}$}
\newtheorem*{conjTfin}{Conjecture T$_{\text{fin}}$}
\newtheorem*{conjNcan}{Conjecture N$_{\text{an}}$}
\newtheorem*{conjNcar}{Conjecture N$_{\text{tb}}$}
\newtheorem{corollary}[theorem]{Corollaire}
\newtheorem{lemma}[theorem]{Lemme}
\newtheorem{proposition}[theorem]{Proposition}
\theoremstyle{definition}
\newtheorem{definition}[theorem]{D\'efinition}

\newtheorem{remark}[theorem]{Remarque}
\numberwithin{equation}{section}
\begin{document}
\title{Fibrations et conjecture de Tate}
\author{Marc Hindry, Am\'{\i}lcar Pacheco et Rania Wazir}
\date{\jour}
\address{Universit\'e Denis Diderot Paris VII\\U.F.R. Math\'ematiques\\  case 7012\\  2 Place Jussieu\\  75251 Paris, 
France}
  \email{hindry@math.jussieu.fr}
\address{Universidade Federal do Rio de Janeiro (Universidade do Brasil)\\ Departamento de Ma\-te\-m\'a\-ti\-ca Pura\\
Rua Guai\-aquil 83, Cachambi, 20785-050 Rio de Janeiro, RJ, Brasil}
\email{amilcar@impa.br}
\address{Universit\`a degli Studi di Torino\\Dipartimento di Matematica\\ via Carlo Alberto, 10\\ 10123 Torino, Italia}
  \email{wazir@dm.unito.it}
\thanks{Am\'{\i}lcar Pacheco a \'et\'e partiellement soutenu par les bourses de recherche CNPq 300896/\hfill 91-3, CNPq 
Edital Universal 470099/2003-8, et par le projet PRONEX 41.96.0830.00.  Rania Wazir a \'et\'e partiellement soutenue par 
le projet EAGER (contract number HPRN-CT-2000-00099) et par le GNSAGA de l'INDAM.  Ce travail a \'et\'e commenc\'e pendant
une visite du deuxi\`eme auteur au premier auteur \`a l'Institut de Math\'ematiques de Jussieu, dans le cadre de l'accord 
Br\'esil-France 69.0014/01-5, les deux auteurs remercient cet accord pour son soutien financier et aussi l'Institut pour 
sa chaleureuse ambiance scientifique.}
\begin{abstract} Nous d\'ecrivons le comportement du rang du groupe de Mordell-Weil de la vari\'et\'e de Picard de la 
fibre g\'en\'erique d'une fibration en termes de contributions locales donn\'ees par des moyennes de traces de Frobenius 
agissant sur les fibres. Les \'enonc\'es fournissent une r\'einterpretation de la conjecture de Tate (pour les diviseurs) 
et g\'en\'eralisent des r\'esultats ant\'erieurs de Nagao, Rosen-Silverman et des auteurs.
\par\bigskip
\noindent\textsc{Abstract.} \textbf{Fibrations and Tate's conjecture.}
We describe the behaviour of the rank of the Mordell-Weil group of
the Picard variety of the
generic fibre of a fibration in terms of local contributions given by
averaging traces of
Frobenius acting on the fibres. The results give a reinterpretation
of Tate's conjecture (for
divisors) and generalises previous results of Nagao, Rosen-Silverman
and the authors.
\end{abstract}
  \maketitle

\section{Introduction et discussion.}
Nous appellerons une \emph{fibration} un morphisme
$f:\cax\rightarrow\cay$ propre et plat et dont
la fibre g\'en\'erique est g\'eom\'etriquement irr\'eductible, entre
deux vari\'et\'es alg\'ebriques, lisses et projectives; on supposera
de plus que $f$, $\cax$ et $\cay$
sont d\'efinis sur un corps de nombres $k$. Par convention les
vari\'et\'es soient connexes dans cet
article. On note
$\ov{k}$ une cl\^oture alg\'ebrique de $k$ et
$G_k:=\gal(\ov{k}/k)$ son groupe de Galois absolu. On notera
$n:=\dim\cax$ et $m:=\dim\cay$. Ainsi
toutes les fibres
$\cax_y:=f^{-1}(y)$ sont de dimension $d:=n-m$ et, pour $y$ hors d'un
diviseur $\Delta$ de $\cay$,
on sait que $\cax_y$ est irr\'eductible et lisse. Si nous avons
besoin de fixer la dimension
relative $d$ de $f$, on appellera cette fibration une
\emph{$d$-fibration}. En particulier la fibre g\'en\'erique, not\'ee
$X$, est une vari\'et\'e
projective lisse de dimension $d$, d\'efinie sur le corps de fonctions
$K:=k(\cay)$. Pour une vari\'et\'e lisse $\cax$, on notera
$\pic(\cax)$ (resp. $\ns(\cax)$) le groupe de Picard
(resp. de N\'eron-Severi) des $\bar{k}$-diviseurs modulo
\'equivalence rationnelle (resp. modulo
\'equivalence alg\'ebrique); le sous-groupe des classes de diviseurs
d\'efinis sur $k$ sera not\'e
$\pic(\cax/k)$ (resp. $\ns(\cax/k)$). Nous noterons $P_X:=\pic^0(X)$
la vari\'et\'e de Picard de
$X$ et $(\tau,B)$ sa $K/k$-trace; c'est-\`a-dire que tout
$K$-morphisme $\phi$ d'une vari\'et\'e ab\'elienne $A$, d\'efinie sur
$k$, vers $P_X$ se factorise
en
$\phi=\tau\circ\phi'$ (voir \cite{lan2}).

Les  Th\'eor\`emes  de Mordell-Weil et Severi, \'etendus par Lang et
N\'eron, permettent d'affirmer que les
groupes
$B(k)$, $\ns(\cax)$, $\ns(\cax/k)$, $\ns(\cay)$, $\ns(\cay/k)$,
$\ns(X)$, $\ns(X$ \linebreak $/K)$,
$P_X(\bar{k}(\cay))/\tau B(\bar{k})$ et $P_X(k(\cay))/\tau B(k)$ sont
des groupes de type fini
(voir \cite{lan1}). Nous allons \'etudier des interpr\'etations (en
bonne partie conjecturales)
des {\it rangs} de ces groupes en termes d'invariants locaux
associ\'es \`a $f$ que nous passons
maintenant \`a d\'efinir.

Soit $\cao_k$ l'anneau des entiers du corps de nombres $k$ et $\frp$
un id\'eal premier de $\cao_k$, on
note
$\bbf_{\frp}$ son corps r\'esiduel et $q_{\frp}:=\#\bbf_{\frp}$. Pour
une vari\'et\'e lisse $\caz/k$ on
notera
$\tilde{\caz}_{\frp}/\bbf_{\frp}$ sa r\'eduction modulo $\frp$ (dans
cette d\'efinition et dans tout ce qui suit
on s'autorise \`a \'eliminer un ensemble fini - qu'on notera $R$ - de
``mauvais" id\'eaux premiers,
qui peut
\^etre
\'elargi d\`es que n\'ecessaire). On notera $\Frob_{\frp}\in G_k$ un
automorphisme de Frobenius associ\'e \`a
$\frp$, agissant sur $\caz\times_k\ov{k}$ \`a travers $\ov{k}$, et
$I_{\frp}\subset G_k$ le groupe
d'inertie associ\'e
\`a
$\frp$ (d\'efini \`a conjugaison pr\`es).

Pour all\'eger les notations, on notera
$$
H^i(\caz):=\,H^i_{\et}(\caz\times_k\ov{k},\bbq_{\ell})
$$
les espaces de cohomologie $\ell$-adique associ\'es sur lesquels
$G_k$ agit, en particulier,
$\Frob_{\frp}$ agit sur
$H^i(\caz)^{I_{\frp}}$, i.e., l'espace invariant par l'action de
$I_{\frp}$, (en toute rigueur on devrait noter
ces espaces $H_{\ell}^i(\caz)^{I_{\frp}}$, mais ils n'interviendront
qu'\`a travers le
poly\-n\^o\-me caract\'eristique de
$\Frob_{\frp}$ et celui-ci est ind\'ependant de $\ell$, au moins pour
$\frp$ hors d'un ensemble fini d'id\'eaux
de
$\cao_k$). Dans tout cet article nous noterons~:
$$
a_{\frp}(\caz)=\tr(\Frob_{\frp}\,|\,H^1(\caz)^{I_{\frp}}) \text{ et }
b_{\frp}(\caz)=\tr(\Frob_{\frp}\,|\,H^2(\caz)^{I_{\frp}}).
$$
Il est classique de d\'efinir les fonctions $L$ associ\'ees \`a
$\caz/k$, que nous ne consid\'ererons que pour
$i=1,2$, par:
\begin{equation}\label{fl}
L_i(\caz/k,s)=\prod_{\frp}\det(1-q_{\frp}^{-s}\Frob_{\frp}\,|\,H^i(\caz)^{I_{\frp}})^{-1}
\end{equation}
Ce produit et la s\'erie de Dirichlet associ\'ee sont convergentes
pour $\Re(s)>1+i/2$.

On note dans ce texte par $F_{\frp}$ le g\'en\'erateur topologique de
$G_{\bbf_\frp}:=\gal(\bar{\bbf}_{\frp}/\bbf_{\frp})$ qu'on appelle le
Frobenius en caract\'eristique $p$. Notons
aussi
$H^i(\tilde{\caz}_{\frp}):=H^i_{\et}(\tilde{\caz}_{\frp}\times_{{\bbf}
_{\frp}}\ov{\bbf}_{\frp},\bbq_{\ell})$.
Observons que les th\'eor\`emes de sp\'ecialisation (ou changement de
base) en cohomologie \'etale
(voir \cite{mil} et
\cite[Appendix C]{ha}) permettent d'identifier le polyn\^ome
caract\'eristique de $F_{\frp}$ agissant
sur
$H^i(\tilde{\caz}_{\frp})$ et celui de $\Frob_{\frp}$ agissant sur
$H^i(\caz)^{I_{\frp}}$; en
particulier
$\tr(F_{\frp}|H^i(\tilde{\caz}_{\frp}))=\tr(\Frob_{\frp}|H^i(\caz)^{I_
{\frp}})$.

Quitte \`a \'elargir $R$, si n\'ecessaire, on peut supposer que la r\'eduction
$\tilde{f}_{\frp}:\tilde{\cax}_{\frp}\to\tilde{\cay}_{\frp}$ de
$f:\cax\to\cay$ soit aussi une fibration
en vari\'et\'es lisses d\'efinie sur $\bbf_{\frp}$. Soit
$\Delta:=\{y\in\cay\,|\,\cax_y=f^{-1}(y)\text{ est singulier}\}$ le
lieu discriminant de $f$. En
rajoutant un nombre fini d'id\'eaux \`a $R$, on peut supposer que la
r\'eduction
$\tilde{\Delta}_{\frp}$ de $\Delta$ modulo $\frp$ co\"{\i}ncide avec
le lieu discriminant de $\tilde{f}_{\frp}$.

Pour tout
$y\in(\tilde{\cay}_{{\frp}}-\tilde{\Delta}_{\frp})(\bbf_{\frp})$, soit
$a_{\frp}(\cax_{\frp,y}):=\tr(F_{\frp}\,|\,H^1(\tilde{\cax}_{\frp,y}))$ et
$b_{\frp}(\cax_{\frp,y}):=\tr(F_{\frp}\,|\,H^2(\tilde{\cax}_{\frp,y}))$. Si
$y\in\tilde{\Delta}_{\frp}(\bbf_{\frp})$, on doit remplacer les
groupes de cohomologie \'etale
$H^i(\tilde{\cax}_{\frp,y})$ par les groupes de cohomologie
$\ell$-adique \`a support propre
$H^i_c(\tilde{\cax}_{\frp,y})$ \linebreak $:=H^i_c(\tilde{\cax}_{\frp,y}\times_{\bbf
_{\frp}}\ov{\bbf}_{\frp},\bbq_{\ell})$,
i.e.,
$a_{\frp}(\cax_{\frp,y}):=\tr(F_{\frp}\,|\,H^1_{c}(\tilde{\cax}_{\frp,y}))$ et
$b_{\frp}(\cax_{\frp,y}):=\tr(F_{\frp}\,|\,H^2_{c}($ \linebreak $\tilde{\cax}_{\frp,y}))$.

Nous introduisons  les
\emph{traces moyennes de Frobenius} d\'efinies par
$$
\frA_{\frp}(\cax)=\frac1{q_{\frp}^{m}}\sum_{y\in\tilde{\cay}_{\frp}(\bbf_{\frp})}a_{\frp}(\tilde{\cax}_{\frp,y}) \text{
et } \frB_{\frp}(\cax)=\frac1{q_{\frp}^{m}}\sum_{y\in\tilde{\cay}_{\frp}(\bbf_{\frp})}b_{\frp}(\tilde{\cax}_{\frp,y}).
$$

  On d\'efinit \'egalement comme
dans \cite{hp} la
\emph{trace moyenne de Frobenius r\'eduite}   par
$$\frA_{\frp}^*(\cax)=\frA_{\frp}(\cax)-a_{\frp}(B).$$

Avant d'\'enoncer le th\'eor\`eme principal de ce texte, rappelons le
contenu de la Conjecture de Tate (pour
les diviseurs).

\begin{conjTa}[Conjecture de Tate]\label{conj1a}\cite[Conjecture 2]{ta1}
La fonction $L_2(\cax/k,s)$ a un p\^ole en $s=2$ d'ordre $\rg(\ns(\cax/k))$.
\end{conjTa}

\begin{remark}
On peut pr\'eciser l'\'enonc\'e comme suit.
\begin{enumerate}
\item Il s'agit d'une version de la conjecture de Tate, la conjecture
g\'en\'erale concerne tous les
cycles alg\'ebriques.

\item On peut se dispenser de l'hypoth\`ese d'un prolongement
m\'eromorphe au voisinage de $s=2$ en interpr\'etant
la phrase ``$L_2(\cax/k,s)$ poss\`ede un p\^ole d'ordre $t$ en $s=2$"
comme signifiant
$$
\lim_{\Re(s)>2, s\to 2}(s-2)^tL_2(\cax/k,s)=\alpha\not= 0.
$$
De m\^eme, si
$f(s)$ est holomorphe sur $\Re(s)>\lambda$ et si
$\lim_{s\to\lambda}(s-\lambda)f(s)=\alpha\not= 0$, on
appellera
$\alpha$ le {\it r\'esidu} de la fonction $f(s)$ en $s=\lambda$ et on \'ecrira
$\res_{s=\lambda}f(s)=\alpha$.
\end{enumerate}
\end{remark}

\begin{conjMan}
  Avec les notations pr\'ecedentes. La fonction
\begin{equation}\label{nagao01}
\sum_{\frp\notin
R}-\frA_{\frp}^*(\cax)\frac{\log(q_{\frp})}{q_{\frp}^s}+\sum_{\frp\not
in
R}\frB_{\frp}(\cax)\frac{\log(q_{\frp})}{q_{\frp}^{s+1}},
\end{equation}
qui est a priori analytique sur $\Re(s)>1$, a un
p\^ole simple en $s=1$, avec r\'esidu  \'egal \`a:
\begin{equation}\label{nagao02}
\rg\left(\frac{P_X(K)}{\tau B(k)}\right)+\rg(\ns(X/K)).
\end{equation}
\end{conjMan}

Nous pouvons maintenant \'enoncer notre th\'eor\`eme principal qui
est une cons\'e\-quen\-ce d'un
r\'esultat inconditionnel de calcul de r\'esidus (cf. \S4, Th\'eor\`eme
\ref{th4ia}).

\begin{theorem}\label{th1a}
Soit $f:\cax\to\cay$ une fibration en vari\'et\'es projectives
lisses, d\'efinies sur un corps de nombres $k$. Deux
des trois affirmations suivantes impliquent la troisi\`eme:
\begin{enumerate}
\item La Conjecture T est vraie pour $\cax/k$.

\item La Conjecture T est vraie pour $\cay/k $.

\item La Conjecture $M_{\text{an}}$ est vraie pour $f:\cax\to\cay$.
\end{enumerate}
\end{theorem}

En fait, conjecturalement chacune des deux fonctions d\'efinies par
les sommes dans (\ref{nagao01}) devrait
se prolonger avec un r\'esidu en $s=1$ \'egal au terme correspondant
dans (\ref{nagao02}). Pour
mettre en perspective l'\'enonc\'e pr\'ec\'edent il convient de
rappeler la Conjecture Analytique
de Nagao G\'en\'eralis\'ee (ainsi  que son analogue arithm\'etique,
la Conjecture
Taub\'erienne de Nagao G\'en\'eralis\'ee) et formuler un analogue de la
Conjecture T pour les fibrations en vari\'et\'es d\'efinies sur un corps de
nombres
$k$.

\begin{conjNcan}[Conjecture Analytique de Nagao G\'en\'eralis\'ee]
\begin{equation}\label{eq1a}
\res_{s=1}\left(\sum_{\frp\notin
R}-\frA_{\frp}^*(\cax)\frac{\log(q_{\frp})}{q_{\frp}^s}\right)=\rg\left(\frac{P_X(K)}{\tau B(k)}\right).
\end{equation}
\end{conjNcan}

\begin{conjNcar}[Conjecture Taub\'erienne de Nagao G\'en\'eralis\'ee]
\begin{equation}\label{eq1b}
\lim_{T\to\infty}\frac1T\left(\sum_{\substack{\frp\notin R\\ q_{\frp}\le
T}}-\frA_{\frp}^*(\cax)\log(q_{\frp})\right)=\rg\left(\frac{P_X(K)}{\tau B(k)}\right).
\end{equation}
\end{conjNcar}

\begin{conjTfin}[Analogue de la Conjecture T pour les fibrations en
vari\'et\'es]
\label{conj1a*}
$$
\res_{s=2}\left(\sum_{\frp\notin
R}\frB_{\frp}(\cax)\frac{\log(q_{\frp})}{q_{\frp}^s}\right)=\rg\left(\ns(X/K)\right).
$$
\end{conjTfin}

On peut bien s\^ur formuler une version \emph{taub\'erienne} de cette
conjecture,
d\'enot\'ee par Conjecture T$_{\text{fin,tb}}$.

\begin{corollary}[du Th\'eor\`eme \ref{th1a}]\label{co1a}
Trois des affirmations suivantes impliquent la quatri\`eme

\begin{enumerate}
\item La Conjecture T est vraie pour $\cax/k$.

\item La Conjecture T est vraie pour $\cay/k$.

\item La Conjecture T$_{\text{fin}}$ est vraie pour $f:\cax\to\cay$

\item La Conjecture N$_{\text{an}}$ est vraie pour $f:\cax\to\cay$.
\end{enumerate}
\end{corollary}

\begin{remark}\label{remnm}

\begin{enumerate}

\item On voit tout de suite que la Conjecture M$_{\text{an}}$ se
r\'eduit \`a la Conjecture
N$_{\text{an}}$ dans le cas de fibrations en courbes, puisque  la Conjecture
T$_{\text{fin}}$ est trivialement vraie. En effet, dans ce cas-l\`a on a $\dim
X=1$, donc $\ns(X)\cong\bbz$.  De plus,  pour tout
$y\in(\tilde{\cay}_{\frp}-\tilde{\Delta}_{\frp})(\bbf_{\frp})$ on a
$b_{\frp}(\tilde{\cax}_{\frp,y})=q_{\frp}$, alors que
$b_p(\tilde{\cax}_{\frp,y})=O(q_{\frp})$ pour $y\in
\tilde{\Delta}_{\frp}(\bbf_{\frp})$.  Donc $\frB_{\frp}(\cax)$ se
d\'ecompose en
somme de deux termes:
$$
{q_{\frp}^m} {\frB}_{\frp}({\cax}) := \sum_{y \in
\tilde{\cay}_{\frp}(\bbf_{\frp})}
b_{\frp}(\tilde{\cax}_{\frp,y})=
q_{\frp}\#(\tilde{\cay}_{\frp}(\bbf_{\frp}))+\sum_{y\in
\tilde{\Delta}_{\frp}(\bbf_{\frp})} (b_{\frp}(\tilde{\cax}_{\frp,y})-q_{\frp})
$$
o\`u (cf. Lemme {\ref{lem2b}}) le premier terme est
$q_{\frp}\left(q_{\frp}^{\dim(\tilde{\cay}_{\frp})}+O\left(q_{\frp}^{\dim(\tilde{\cay}_{\frp})-1/2 }\right) \right)$
et le second terme est $O \left(q_{\frp}^{1+\dim(\tilde{\Delta}_{\frp})}\right)$.  Il suit que ${\frB}_{\frp}({\cax}) =
q_{\frp} + O\left(q_{\frp}^{1/2}\right)$, et en particulier,
$$
\res_{s=2}\left(\sum_{\frp\notin
R}\frB_{\frp}(\cax)\frac{\log(q_{\frp})}{q_{\frp}^s}\right)=1.
$$

\item La Conjecture T pour $\cax/k$ peut s'exprimer en disant que la
d\'eriv\'ee logarithmique de la
fonction
$L_2(\cax/k,s)$ poss\`ede un p\^ole simple en $s=2$ de r\'esidu
$-\rg\left(\ns(\cax/k)\right)$ ou encore (voir
le paragraphe 4 pour le d\'etail de ces calculs classiques) que~:
\begin{equation}\label{tan}
\res_{s=2}\left(\sum_{\frp\notin
R}b_{\frp}(\cax)\frac{\log(q_{\frp})}{q_{\frp}^s}\right)=\rg(\ns(\cax/k)).
\end{equation}

\item Ainsi, dans la plupart des cas o\`u l'on sait d\'emontrer
l'existence d'un prolongement analytique de
$L_2(\cax/k,s)$
\`a la droite $\Re(s)=2$, on sait \'egalement d\'emontrer que la
fonction ne s'annule pas sur cette droite.
Cette propri\'et\'e est importante car elle permet d'appliquer un th\'eor\`eme
Taub\'erien
\cite[Chapter XV]{la} \`a la d\'eriv\'ee logarithmique de
$L_2(\cax/k,s)$ et de conclure que
\begin{equation}\label{tar}
\lim_{T\to\infty}\frac1T\left(\sum_{\substack{\frp\notin R\\ q_{\frp}\le
T}}b_{\frp}(\cax)\frac{\log(q_{\frp})} {q_{\frp}}\right)=\rg(\ns(\cax/k)).
\end{equation}

Si cette derni\`ere \'egalit\'e est vraie, nous dirons  que la
\emph{Conjecture
Taub\'eri\-en\-ne de Tate}  est v\'erifi\'ee, d\'enot\'ee Conjecture
T$_{\text{tb}}$.

\item Sous l'hypoth\`ese additionnelle que la Conjecture
T$_{\text{tb}}$ soit vraie pour $\cax/k$ et $\cay/k$, et que
la Conjecture T$_{\text{fin,tb}}$ soit vraie pour $f:\cax\to\cay$, on
en conclut que la Conjecture
N$_{\text{tb}}$ est aussi vraie pour $f:\cax\to\cay$.
\end{enumerate}
\end{remark}

Comme nous l'avons indiqu\'e, la Conjecture T$_{\text{fin}}$ est
v\'erifi\'ee (et essentiellement triviale) dans le cas o\`u $X$
est une courbe; dans ce cas sa vari\'et\'e de Picard s'identifie \`a
sa vari\'et\'e Jacobienne que
nous d\'enotons par
$J_X$. On obtient alors~:

\begin{corollary}\label{co1a2}
Soit $f:\cax\rightarrow\cay$ une fibration en courbes. Supposons que
la Conjecture T soit vraie pour $\cax/k$
et
$\cay/k$. Alors

\begin{equation}\label{eq1a2}
\res_{s=1}\left(\sum_{\frp\notin
R}-\frA_{\frp}^*(\cax)\frac{\log(q_{\frp})}{q_{\frp}^s}\right)=\rg\left(\frac{J_X(K)}{\tau B(k)}\right).
\end{equation}

Sous l'hypoth\`ese additionnelle que les fonctions $L_2(\cax/k,s)$ et
$L_2(\cay/k,s)$ se prolongent sur la
droite
  $\Re(s)=2$ sans z\'eros, on obtient \'egalement la
  version taub\'eri\-en\-ne suivante~:
\begin{equation}\label{eq1b2}
\lim_{T\to\infty}\frac1T\left(\sum_{\substack{\frp\notin R\\ q_{\frp}\le
T}}-\frA_{\frp}^*(\cax)\log(q_{\frp})\right)=\rg\left(\frac{J_X(K)}{\tau B(k)}\right).
\end{equation}
\end{corollary}

\begin{remark}
Si l
'on prend une vari\'et\'e ab\'elienne ``constante", disons $B$
d\'efinie sur un corps de nombre $k$, le pendant de
la formule (\ref{eq1a2}) est
$$
\res_{s=1}\left(\sum_{\frp\notin
R}-a_{\frp}(B)\frac{\log(q_{\frp})}{q_{\frp}^s}\right)\overset{?}=\rg(B(k))
$$
c'est-\`a-dire la Conjecture de Birch \& Swinnerton-Dyer pour $B/k$
(voir par exemple Tate \cite{ta2}). On
peut expliciter cette analogie en associant \`a la vari\'et\'e $X/K$
(ou plus pr\'ecis\'ement \`a
la fibration
$f:\cax\rightarrow\cay$) les fonctions $L$ suivantes~:
\begin{equation}\label{flK}\begin{aligned}
L_i(X/K,s)&=\prod_{\frp}\left(\prod_{y\in(\tilde{\cay}_{\frp}-\tilde{\Delta}_{\frp})(\bbf_{\frp})}
\det(1-F_{\frp}q_{\frp}^{-s}
\,|\,H^i(\tilde{\cax}_{\frp,y}))^{-1}\right.\\
&\times\left.\prod_{y\in\tilde{\Delta}_{\frp}(\bbf_{\frp})}\det(1-F_{\frp}q_{\frp}^{-s}
\,|\,H^i_{c}(\tilde{\cax}_{\frp,y}))^{-1}\right),
\end{aligned}\end{equation}

Un calcul simple montre alors que la formule analogue de Birch \&
Swinnerton-Dyer pour les fibrations en
vari\'et\'es:
$$
\ord_{s=m+1}L_1(X/K,s)\overset{?}=\rg(P_X(K))
$$
  est \'equivalente \`a la formule
\begin{equation}\label{eqbsdfin}
\res_{s=1}\left(\sum_{\frp\notin
R}-\frA_{\frp}(\cax)\frac{\log(q_{\frp})}{q_{\frp}^s}\right)\overset{?
}=\rg(P_X(K))
\end{equation}
elle-m\^eme cons\'equence de la Conjecture de Birch \&
Swinnerton-Dyer pour $B/k$ et de la formule (\ref{eq1b}).
On appelle la validit\'e de la formule (\ref{eqbsdfin}) la
\emph{Conjecture BSD$_{\text{fin}}$}
pour la fibration
$f:\cax\to\cay$ (cf. \cite[Conjecture of B + SD]{ta1}).

De m\^eme on peut ais\'ement montrer que la Conjecture
T$_{\text{fin}}$ pour une fibration $f:\cax\to\cay$ peut
se reformuler en disant que
$$
\ord_{s=m+2}(L_2(X/K,s))=\rg(\ns(X/K))
$$
(cf. \cite[Conjecture 2]{ta1}). Avec ces notations, la Conjecture
M$_{\text{an}}$, et donc la
troisi\`eme affirmation du Th\'eor\`eme \ref{th1a}, se traduit en disant que
$$
\ord_{s=1}\left(\frac{L_1(X/K,m+s)}{L_1(B/k,s)L_2(X/K,m+1+s)}\right)
$$
est \'egal \`a la quantit\'e attendue~:
$$
\rg\left(\frac{P_X(K)}{\tau B(k)}\right)+\rg\left(\ns(X/K)\right).
$$
\end{remark}

\begin{remark}
Il peut sembler plus naturel de d\'efinir
$$
\frA'_{\frp}(\cax)=\frac1{\#\tilde{\cay}_{\frp}(\bbf_{\frp})}\sum_{y\in
\tilde{\cay}_{\frp}(\bbf_{\frp})}a_{\frp}(\tilde{\cax}_{\frp,y})
$$
et $\frA^{'*}_{\frp}(\cax)=\frA'_{\frp}(\cax)-a_{\frp}(B)$. En
observant (cf. \S3) que l'on dispose de
l'estimation
$\#\tilde{\cay}_{\frp}(\bbf_{\frp})=q_{\frp}^{m}+O(q_{\frp}^{m-1/2})$,
on peut montrer, comme dans
\cite[Remarque 2.3]{hp} que cela ne modifierait pas les \'enonc\'es en vue.
\end{remark}

\begin{remark}
L'\'egalit\'e
  (\ref{eq1b}) ou (\ref{eq1b2}) est une g\'en\'eralisation de la
\emph{Conjecture
Taub\'erienne de Nagao} \cite{nag1,nag2} qui concernait les
fibrations en courbes
elliptiques au dessus de la droite projective et $k=\bbq$ (voir aussi
\cite[Nagao's
Conjecture 1.1$^\prime$]{rosi}). L'\'egalit\'e (\ref{eq1a}) ou
(\ref{eq1a2}) est
une g\'en\'eralisation de la \emph{Conjecture Analytique de Nagao}
telle qu'elle a
\'et\'e formul\'ee  par Rosen et Silverman (cf. \cite[Nagao's Conjecture
1.1]{rosi}) dans le cas o\`u $n=2$ (donc la Conjecture T est trivialement
v\'erifi\'ee pour $\cay$), le genre arithm\'etique des fibres est
\'egal \`a 1 et
$f$ admet une section. De plus, le sch\'ema de leur preuve a servi de
mod\`ele pour les travaux
ult\'erieurs des auteurs qui ont respectivement trait\'e le cas d'une
vari\'et\'e $\cax$ de
dimension 3 fibr\'ee en courbes elliptiques au-dessus d'une surface
\cite[Theorem 1.1]{wa}, d'une
surface fibr\'ee en courbes de genre quelconque
\cite[Th\'eor\`eme 1.3]{hp}.  On peut \'egalement citer le travail de
S. Wong \cite[Theorem 5]{sw}
qui a prouv\'e un analogue  de \cite[Th\'eor\`eme 1.3]{hp},  sous les
hypoth\`eses additionnelles
que $B=0$, et les multiplicit\'es des composantes des fibres
singuli\`eres sont premi\`eres entre
elles. Dans le cadre des surfaces elliptiques, Silverman
\cite[Theorem 6]{si} calcule une borne sup\'erieure pour
$|\frA_{\frp}(\cax)|$, qui d\'epend du
degr\'e du conducteur de $X/K$, du genre de $C/k$ et d'un terme
d'erreur $O(q_{\frp}^{-1/2})$.
\end{remark}

Le texte est organis\'e comme suit. Le \S2 contient les
pr\'eliminaires g\'eom\'etriques aboutissant \`a la
``Formule de Shioda-Tate" (Proposition \ref{prop2b}) dans notre
contexte; le \S3 d\'ecrit
l'estimation du nombre de points sur les corps finis de vari\'et\'es
lisses ou singuli\`eres, le
point clef \'etant la Proposition
\ref{sfc}; le
\S4 contient la preuve du th\'eor\`eme dont d\'erive imm\'ediatement
le Th\'eor\`eme \ref{th1a}.
Dans le \S5 on utilise un r\'esultat g\'eom\'etrique (cf. Proposition
\ref{prop5a}) pour prouver
l'\'equivalence entre la Conjecture T, la Conjecture M$_{\text{an}}$
pour les fibrations en
vari\'et\'es et la Conjecture N$_{\text{an}}$ pour les fibrations en
courbes. Finalement, dans le
\S6 on exploite le rapport entre les Conjectures N$_{\text{an}}$ et
T$_{\text{fin}}$ pour les
fibrations en vari\'et\'es montrant qu'elles sont en effet
\'equivalentes. On prouve aussi que la Conjecture BSD$_{\text{fin}}$
pour les fibrations en vari\'et\'es implique
la Conjecture T$_{\text{fin}}$ pour les fibrations en vari\'et\'es.




\section{Outils g\'eom\'etriques}

Nous prouvons ici une formule (Proposition \ref{prop2b}) analogue \`a
celle de Shioda-Tate  pour les
surfaces elliptiques ayant une section g\'en\'eralisant \cite{rosi},
\cite{wa}, \cite{hp}. Nous
utilisons la th\'eorie d'intersection et prenons comme r\'ef\'erence
Fulton \cite{fu} ou
Hartshorne \cite[Appendix A]{ha}. Un point d\'elicat est de prouver
que l'accouplement  d\'efini
avant le  Lemme \ref{lem2c} est d\'efini n\'egatif.

Soit $f^*:\Div(\cay)\to\Div(\cax)$ l'application pullback des
diviseurs et $f^*:\pic(\cay)\to\pic(\cax)$
l'application d\'eduite pour les classes. Soit $\ch^i(\cax)$ le
groupe de Chow des cycles de
codimension $i$ (donc
$\ch^1(\cax)=\pic(\cax)$). Un \'el\'ement de $\pic(\cay)$ est {\it
positif} (resp. {\it strictement positif}) s'il
est la classe d'un diviseur effectif (resp. effectif non nul).
Choisissons une section hyperplane
de $\cax$ et notons
$h$ sa classe dans $\ch^1(\cax)=\pic(\cax)$. On d\'efinit un accouplement
$\langle\cdot,\cdot\rangle$ par~:

\begin{equation*}\begin{matrix}
\pic(\cax)\times\pic(\cax)&\rightarrow&\ch^2(\cax)&\rightarrow&\ch^{n-
m+1}(\cax)
&\overset{f_*}\longrightarrow&\pic(\cay)\\
(\xi,\zeta)&\mapsto&\xi.\zeta&\mapsto&\xi.\zeta.h^{n-m-1}&\mapsto&\langle\xi,\zeta\rangle:=\\
&&&&&&f_*(\xi.\zeta.h^{n-m-1})
\end{matrix}\end{equation*}

On dira qu'un diviseur de $\cax$ est {\it vertical} si ses
composantes sont contenues dans $f^{-1}(Z)$ pour $Z$
un diviseur de $\cay$; en d'autres termes~: un diviseur
irr\'eductible $D$ de $\cax$ est vertical
si $f(D)\not=\cay$ et un diviseur est vertical s'il est somme de
diviseurs verticaux
irr\'eductibles (bien entendu cette notion est relative
\`a la fibration $f$). Un diviseur irr\'eductible $D$ de $\cax$ sera
dit {\it horizontal} si
$f(D)=\cay$, un diviseur horizontal est une somme de tels diviseurs.
Tout diviseur ou classe de
diviseur (dans
$\pic(\cax)$ ou $\ns(\cax)$) se d\'ecompose en somme d'un diviseur
vertical et d'un autre
horizontal, en particulier nous pouvons \'ecrire (avec des notations
\'evidentes) une
d\'ecomposition de $\bbz[G_k]$-modules
\begin{equation}\label{eqns}
\ns(\cax)=\ns_{\rm hor}(\cax)\oplus \ns_{\rm ver}(\cax).
\end{equation}

Pour le cas o\`u les fibres de $f$ sont des courbes de  genre
arithm\'etique  \'egal \`a 1 et $f$
admet une section, le  r\'esultat suivant est donn\'e dans la th\`ese
du troisi\`eme auteur, o\`u
l'accouplement d\'efini pr\'ec\'edemment est introduit
\cite[Proposition 2.2]{wa}.

\begin{lemma}\label{lem2c}
L'accouplement est n\'egatif sur les diviseurs verticaux, i.e.,  si
$V$ est vertical, on a $\langle V,V\rangle\leq
0$. De plus $\langle V,V\rangle=0$ si et seulement si il existe
$D\in\pic(\cay)\otimes\bbq$ tel que
$V=f^*(D)$. En particulier, si $V$ est vertical et num\'eriquement
(ou alg\'ebriquement) nul sur
$\cax$ alors il appartient \`a
$f^*(\pic(\cay)\otimes\bbq)$.
\end{lemma}

\begin{proof}[\Dem]
On prouve en premier que $\langle V,V\rangle\le 0$. \'Ecrivons
$V=\sum_i V_i$ avec
$V_i\subset f^{-1}(Z_i)$ et $Z_i$
irr\'eductible. Donc, $\langle V_i,V_j\rangle=0$ pour $i\ne j$ et
$\langle V,V\rangle=\sum_i\langle V_i,V_i\rangle$. Il
suffit donc de prouver l'affirmation pour chaque $V_i$. Par cons\'equent on
peut supposer que $V\subset f^{-1}(G)$ pour
$G\in\Div(\cay)$ irr\'eductible. D\'ecomposons $\cag=f^*(G)=\sum n_i\gamma_i$
avec $n_i\ge 1$ pour chaque $i$, et
$\gamma_i$ irr\'eductible dans $\cax$. En plus, comme $V\subset f^{-1}(G)$,
on obtient $V=\sum_i a_i\gamma_i$ qu'on
re\'ecrit comme $\sum_i\frac{a_i}{n_i}n_i\gamma_i$.
Soit $V'=\sum_i(\frac{a_i}{n_i})^2n_i\gamma_i$. Comme $V'$ est
vertical et $\cag=f^*(G)$, on a $\langle V',\cag\rangle=0$. Donc,
\begin{align*}
-2\langle V,V\rangle&=\langle V',\cag\rangle-2\langle
V,V\rangle+\langle\cag,V'\rangle\\
&=\sum_{i,j}\left(\frac{a_i}{n_i}-\frac{a_j}{n_j}\right)^2\langle
n_i\gamma_i,n_j\gamma_j\rangle
\end{align*}
et a fortiori
$$\langle V,V\rangle=-\frac12\sum_{i\ne
j}\left(\frac{a_i}{n_i}-\frac{a_j}{n_j}\right)^2
\langle n_i\gamma_i,n_j\gamma_j\rangle.$$ Par la d\'efinition de
l'intersection avec un diviseur \cite[Chapter 2,
2.3, p. 33]{fu} $\gamma_i.\gamma_j\ge 0$ pour $i\ne j$ et donc, comme
$h$ est la classe d'un
diviseur tr\`es ample
$\gamma_i.\gamma_j.h^{n-m-1}\ge 0$ (cf. \cite[Lemma 12.1, p.
211]{fu}). Par la d\'efinition de l'image directe
on obtient $\langle\gamma_i,\gamma_j\rangle\ge 0$ pour $i\ne j$. On
en d\'eduit $\langle
n_i\gamma_i,n_j\gamma_j\rangle\ge 0$, d'o\`u $\langle V,V\rangle\le 0$.

Supposons maintenant que $\langle V,V\rangle=0$. On a alors
$\frac{a_i}{n_i}=\frac{a_j}{n_j}$ pour chaque $i\ne j$
tels que $\langle\gamma_i,\gamma_j\rangle>0$.  Nous
affirmons que, pour chaque $i\ne j$ il existe une cha\^{\i}ne de
composantes irr\'eductibles
$\gamma_i=\gamma_{i_0},\cdots,\gamma_{i_r}=\gamma_j$ telles que
$\langle \gamma_{i_k}, \gamma_{i_{k+1}} \rangle$ est $>0$. On en
d\'eduira donc que $\frac{a_i}{n_i}=\frac{a_j}{n_j}$ pour tout $i\ne
j$ et ainsi
que
$V=\sum_ia_i\gamma_i=\frac{a_1}{n_1}\sum_in_i\gamma_i=\frac{a_1}{n_1}f
^*(G)$ comme annonc\'e.
Pour v\'erifier
l'affirmation, consid\'erons $\eta_G$ le point g\'en\'erique de $G$
(qui est une vari\'et\'e irr\'eductible) et
la fibre $\cax_{\eta_G}=f^{-1}(\eta_G)$. Celle-ci est une union de
composantes irr\'eductibles
$X_1\cup\dots \cup X_r$ et on voit ais\'ement que les $X_{\nu}$
correspondent aux $\gamma_{\nu}$ :
en fait $\gamma_{\nu}$ est la fermeture de Zariski de $X_{\nu}$ et
$X_{\nu}=\gamma_{\nu}\cap\cax_{\eta_G}$. En outre, comme $f$ est plat,
$\cax_{\eta_G}$ est connexe, et en particulier pour chaque $\nu \ne
\mu$ il existe une cha\^{\i}ne
de composantes irr\'eductibles
$X_{\nu}=X_{\nu_0},\cdots, X_{\nu_r}= X_{\mu}$ telles que $X_{\nu_k}
\cap X_{\nu_{k+1}} \ne \emptyset$.  Donc $\gamma_{\nu_k} \cap
\gamma_{\nu_{k+1}} \ne \emptyset$.  Soient alors $H_1,\dots,H_{n-m-1}$
des repr\'esentants de $h$ dans $\Div(\cax)$ en position suffisamment
g\'en\'erale. Le morphisme
$f:\gamma_{\nu_k}\cap\gamma_{\nu_{k+1}}\cap H_1\cap \dots\cap
H_{n-m-1}\rightarrow \cay$ (\`a valeurs dans $G$)
est g\'en\'eriquement fini car au-dessus de $\eta_G$ on trouve
$X_{\nu_k}\cap X_{\nu_{k+1}}\cap
H_1\cap \dots\cap H_{n-m-1}$ qui est fini et non vide (pour des
raisons de dimension). Ainsi, on
trouve que
$f_*(\gamma_{\nu_k}\cdot \gamma_{\nu_{k+1}}\cdot h^{n-m-1})>0$, ce
qu'il fallait d\'emontrer.
\end{proof}

L'\'etude de la $K/k$-trace $B$ de $P_X$ est similaire \`a
\cite[\S3]{hp}. L'application pullback induit
un homomorphisme que nous notons $f^*:\pic^0(\cay)\to\pic^0(\cax)$.
L'inclusion de la fibre
g\'en\'erique
$\iota:X\rightarrow\cax$ induit un homomorphisme ``restriction \`a la
fibre g\'en\'erique"
$\iota^*:\pic^0(\cax)\to\pic^0(X)$ qui, par la d\'efinition de la
$K/k$-trace, se factorise par un
homomorphisme
$b:\pic^0(\cax)\to B$ tel que $\tau\circ b=\iota^*$.

\begin{proposition}(Cf. \cite[Proposition 3.4]{hp})\label{prop2c}
Le groupe $\ker(f^*)$ est fini et est m\^eme trivial si $f$ admet une
section g\'en\'erique. La suite  de
vari\'et\'es ab\'eliennes
\begin{equation}\label{eq2f}
0\to\pic^0(\cay)\to\pic^0(\cax)\to B\to0.
\end{equation}
est exacte si $f$ poss\`ede une section g\'en\'erique et exacte \`a
des groupes finis pr\`es en g\'en\'eral (ou
encore apr\`es tensorisation par $\bbq$). En particulier,
$H^1(\cax)\cong H^1(\cay)\oplus H^1(B)$.
A fortiori,  pour presque tout id\'eal premier $\frp$ de $\cao_k$, on a
$a_{\frp}(\cax)=a_{\frp}(\cay)+a_{\frp}(B)$.
\end{proposition}

\begin{proof}[\Dem]
La preuve est essentiellement la m\^eme que dans \cite[Proposition
3.4]{hp}, il s'agit seulement d'ajouter
que d'apr\`es le Lemme \ref{lem2c}, si $V$ est un diviseur vertical
alg\'ebriquement \'equivalent
\`a $0$, c'est une somme de fibres, c'est-\`a-dire qu'il est de la
forme $f^*(D')$ avec
$D'\in\pic^0(\cay)\otimes\bbq$ ou encore il existe
$m\in\bbn^*$ tel que $mV\in f^*\pic^0(\cay)$. L'entier $m$ est
born\'e par les multiplicit\'es des fibres
et m\^eme on peut prendre $m=1$, s'il y a une section g\'en\'erique
(en effet cette derni\`ere est alors
d\'efinie sur le compl\'ementaire d'un ferm\'e de codimension au
moins deux). L'isomorphisme
$H^1(\cax)\cong H^1(\cay)\oplus H^1(B)$ suit de
\cite[Corollary 4.19, p. 131]{mil}.
\end{proof}

Nous passons maintenant \`a (la g\'en\'eralisation de) la formule de
Shioda-Tate. Soit
$\imath^*:\Div(\cax)\to\Div(X)$ l'application pullback des diviseurs
induite par
$\imath:X\to\cax$. Soit maintenant $Z_1,\dots,Z_r$ des diviseurs sur
$X/\bar{k}(\cay)$ dont les classes forment une base de
$\ns(X/\bar{k}(\cay))\otimes\bbq$ et soient
$\bar{Z}_1,\dots,\bar{Z}_r$ les diviseurs sur $\cax$ obtenus en
prenant l'adh\'erence de Zariski de
$Z_1,\cdots,Z_r$ (resp.), de sorte que $\imath^*(\bar{Z}_i)=Z_i$. Si
$D\in\Div(\cax)$ on a donc
$\imath^*(\cl(D))$ alg\'ebriquement
\'equivalent \`a $n_1Z_1+\dots+n_rZ_r$ (avec $n_i\in\bbq$). On
d\'efinit alors $\psi:\pic(\cax)\to
P_X(\bar{k}(\cay))\otimes\bbq$ par
$$\psi(\cl(D))=\imath^*\left(D-(n_1\bar{Z}_1+\dots+n_r\bar{Z}_r)\right).$$
On obtient alors d'apr\`es le Lemme \ref{lem2d} ci-dessous une suite exacte
\begin{equation}\label{eq2g}
0\to\ker(\psi)\to\pic(\cax)\otimes\bbq\overset{\psi}\longrightarrow
P_X(\bar{k}(\cay))\otimes\bbq\to0.
\end{equation}
dont nous allons
d\'ecrire le noyau.

\begin{definition}
Soit $\tilde{\cas}$ (respectivement $\cas$) le sous-groupe de
$\pic(\cax)$ (respectivement de $\ns(\cax)$)
engendr\'e par les classes  $\bar{Z}_1,\dots,\bar{Z}_r$ et des
composantes  des images inverses de
diviseurs de $\cay$ par $f$. On notera respectivement
$\tilde{\cas}_{\bbq}$ et $\cas_{\bbq}$ les espaces vectoriels obtenus
par tensorisation par $\bbq$.
\end{definition}

\begin{lemma}\label{lem2d} (Cf. \cite[Lemme 3.8]{hp})
L'application $\psi$ est surjective et son noyau $\ker(\psi)$ est
\'egal \`a $\tilde{\cas}$. Une base
de
${\cas}_{\bbq}$ est fournie par les classes de
$\bar{Z}_1,\dots,\bar{Z}_r$, l'image r\'eciproque par $f$
d'une base de $\ns(\cay)$ et  par les composantes, sauf une, de
chaque image r\'eciproque
$f^{-1}(G)$, pour
$G$ diviseur irr\'eductible de
$\cay$ tel que le nombre des composantes irr\'eductibles de
$f^{-1}(G)$ soit $\ge2$. Soient
$\Delta_1,\dots,\Delta_s$ les diviseurs irr\'eductibles de $\cay$
tels que le nombre $m_i$ des
composantes de
$f^{-1}(\Delta_i)$ soit $\geq2$. On a en particulier,
$$\rg({\cas})=\sum_{i=1}^s(m_i-1)+\rg(\ns(\cay))+\rg(\ns(X/\bar{k}(\cay)).$$
\end{lemma}

\begin{proof}[\Dem]
Soit $z\in P_X$ tel que $z$ repr\'esente la classe $\cl(E)$ d'un
diviseur $E\in\Div(X)$. Soit $\ov{E}$ la
cl\^oture de Zariski de $E$ dans $\cax$. C'est un diviseur de Weil
(donc de Cartier, puisque $\cax$
est lisse) et il est imm\'ediat qu'on a $\psi(\cl(\ov{E}))=z$. Ainsi
$\psi$ est bien surjective.

Il est clair que $\tilde{\cas}\subset\ker(\psi)$. Supposons que
$\cl(D)\in\ker(\psi)$. Alors, il existe une
fonction
$x\in(\bar{k}(\cay))(X)$ telle que
$\divi(x)=\imath^*(D-(n_1\bar{Z}_1+\dots+n_r\bar{Z}_r))$. Soit
$\tilde{x}$ une fonction rationnelle sur $\cax$ telle que
$\tilde{x}_{|X}=x$. En particulier,
$\imath^*(\divi(\tilde{x}))=\divi(x)$. Il est clair que, au niveau
des diviseurs, $\ker(\imath^*)$
est engendr\'e par  les classes des diviseurs irr\'eductibles
verticaux pour la fibration $f$. Donc
$\divi(\tilde{x})-(D-(n_1\bar{Z}_1+\dots+n_r\bar{Z}_r))$ est une
somme des composantes de ces
fibres. On en tire bien que $\cl(D)\in\tilde{\cas}$ et, par suite
$\ker(\psi)=\tilde{\cas}$.

Pour la derni\`ere affirmation, il est clair que les classes de
l'\'enonc\'e forment un syst\`eme g\'en\'erateur
car toutes les fibres sont alg\'ebriquement \'equivalentes. Le fait
qu'elles soient
ind\'ependantes se v\'erifie en utilisant la th\'eorie d'intersection
et le Lemme \ref{lem2c}. La
formule donnant le rang de
$\cas$ est alors imm\'ediate.
\end{proof}

\begin{definition}\label{defvert}
Posons $\caf:=\ns_{\rm ver}(\cax)/f^*(\ns(\cay))$ et
$\caf_{\bbq}:=\caf\otimes\bbq$.
\end{definition}

Remarquons que, avec les notations du lemme pr\'ec\'edent, on a
$\rg(\caf)=\sum_{i}(m_i$ \linebreak $-1)$ et
que, apr\`es avoir tensoris\'e par $\bbq$ nous avons un isomorphisme
de $\bbq[G_k]$-modules
$$\ns_{\rm ver}(\cax)\otimes\bbq\cong
\caf_{\bbq}\oplus\left(\ns(\cay)\otimes\bbq\right).$$
Nous pouvons maintenant \'enoncer (comparer avec \cite[Formula
1.4]{rosi}, \cite[Proposition
3.9]{hp}, \cite[Theorem 4.5]{wa})~:

\begin{proposition}[Formule de Shioda-Tate]\label{prop2b}
Il existe un isomorphisme de $\bbq[G_k]$-modules
\begin{align*}
&\ns(\cax)\otimes\bbq\cong\left(\left(\frac{P_X(\bar{k}(\cay))}
{\tau B(\bar{k})}\right)\otimes\bbq\right)\oplus\cas_{\bbq}\\
&\cong\left(\left(\frac{P_X(\bar{k}(\cay))} {\tau
B(\bar{k})}\right)\otimes\bbq\right)\oplus\left(\ns(\cay)\otimes\bbq\right)\oplus(\caf\otimes\bbq)
\oplus\left(\ns(X/\bar{k}(\cay))\otimes\bbq\right)
\end{align*}

  En particulier, on a
\begin{equation}\begin{aligned}\label{eq2h}
&\rg(\ns(\cax/k))=\rg\left(\frac{P_X(k(\cay))}{\tau
B(k)}\right)+\rg(\cas^{G_k})\\
=&\rg\left(\frac{P_X(k(\cay))}{\tau B(k)}\right)
+\rg(\ns(\cay/k))+\rg(\caf^{G_k})+\rg(\ns(X/k(\cay))).
\end{aligned}
\end{equation}
\end{proposition}

\begin{proof}[\Dem]
Comme dans \cite[Proposition 3.9]{hp} le r\'esultat suit de
(\ref{eq2g}), (\ref{eq2f}) et du Lemme \ref{lem2d}.
La derni\`ere formule s'obtient en prenant les $G_k$-invariants et en
observant  que
$\ns(\cax/k)_{\bbq}=\ns(\cax)_{\bbq}^{G_k}$ et
$\left(P_X(\bar{k}(\cay))/\tau
B(\bar{k})\right)^{G_k}_{\bbq}=\left(P_X(k(\cay))/\tau
B(k)\right)_{\bbq}$.
\end{proof}




\section{D\'ecompte des points sur les corps finis}

Tout comme dans Rosen-Silverman \cite{rosi} et dans les travaux
ult\'erieurs des auteurs \cite{hp} et
\cite{wa}, l'id\'ee de la preuve du th\'eor\`eme principal est de
compter le nombre des points
rationnels de
$\tilde{\cax}_{\frp}$ sur $\bbf_{\frp}$ de deux mani\`eres. La
premi\`ere consiste \`a employer la
formule de Lefschetz pour
$\tilde{\cax}_{\frp}$, la deuxi\`eme \`a compter le nombre des points
dans chaque fibre $\tilde{\cax}_{\frp,y}$
pour
$y\in \tilde{\cay}_{\frp}(\bbf_{\frp})$ et calculer la somme.

Commen\c{c}ons par expliciter l'application de la formule de Lefschetz.

\begin{lemma}\label{lef}
Soit $\cax$ une vari\'et\'e projective lisse de dimension $n$,
d\'efinie sur le corps de nombre $k$, soit ${\frp}$
un id\'eal premier de $\cao_k$ de bonne r\'eduction et soit
$\tilde{\cax}_{\frp}$ la r\'eduction de
$\cax$ en
${\frp}$ d\'efinie sur le corps fini $\bbf_{\frp}$ de cardinal
$q_{\frp}$, alors~:
$$
\#\tilde{\cax}_{\frp}(\bbf_{\frp})=q_{\frp}^n-a_{\frp}(\cax)q_{\frp}^{
n-1}+b_{\frp}(\cax)q_{\frp}^{n-2}
+O(q_{\frp}^{n-3/2}).
$$
\end{lemma}

\begin{proof}[\Dem]
Il s'agit d'une application directe de la \emph{Formule de Lefschetz}
qui, jointe \`a l'observation
\mbox{$\tr(F_{\frp}|H^i(\tilde{\cax}_{\frp}))
=\tr(\Frob_{\frp}|H^i(\cax)^{I_{\frp}})$} de
l'introduction, indique que
$$
\#\tilde{\cax}_{\frp}(\bbf_{\frp})=
\sum_{i=0}^{2n}(-1)^i\tr(\Frob_{\frp}|H^i(\cax)^{I_{\frp}})
$$
en tenant compte d'une part de l'Hypoth\`ese de Riemann (Th\'eor\`eme
de Deligne, \cite[Th\'eor\`eme 1.6]{del74})
qui entra\^{\i}ne que
$\tr(\Frob_{\frp}|H^{i}(\cax)^{I_{\frp}})=O(q_{\frp}^{i/2})$ et,
d'autre part,
de la dualit\'e de Poincar\'e qui donne
$\tr(\Frob_{\frp}|H^{2n-1}(\cax)^{I_{\frp}})=
q_{\frp}^{n-1}\tr(\Frob_{\frp}|$ \linebreak $H^1(\cax)^{I_{\frp}})=q_{\frp}^{n-1}a_{\frp}(\cax)$ et
$\tr(\Frob_{\frp}|H^{2n-2}(\cax)^{I_{\frp}})=
q_{\frp}^{n-2}\tr(\Frob_{\frp}|H^2(\cax)^{I_{\frp}})=q_{\frp}^{n-2}b_{
\frp}(\cax)$. Enfin bien s\^ur on a
l'\'egalit\'e
$\tr(\Frob_{\frp}|H^{2n}(\cax)^{I_{\frp}})=q_{\frp}^{n}$.
\end{proof}

Cependant nous aurons besoin de compter les points sur des fibres non
ir\-r\'e\-du\-ctibles ou singuli\`eres;
une premi\`ere \'etape est donn\'ee par les estimations de Lang-Weil.
Nous appelons ici {\it ensemble
alg\'ebrique} sur $\bbf_{\frp}$ un sch\'ema projectif r\'eduit sur
$\bbf_{\frp}$.

\begin{lemma}\label{lem2b}
Soit $\delta$ un ensemble alg\'ebrique de dimension pure $m$ d\'efini
sur $\bbf_{\frp}$ et soit
$\delta=\delta_1\cup\cdots\cup\delta_r$ sa d\'ecomposition en
composantes irr\'eductibles absolues
telles que
$\delta_i$ soit d\'efinie sur $\bbf_{\frp}$ pour $1\le i\le s$ (et
seulement pour ces valeurs). On a
\begin{equation}\label{eq2c}
\#\delta_i(\bbf_{\frp})=\left\{\begin{aligned}
q_{\frp}^m+O(q_{\frp}^{m-1/2}),&\quad\text{ si }1\le i\le s\\
O(q_{\frp}^{m-1}),&\quad\text{ si }s<i\le r,
\end{aligned}\right.\end{equation}
o\`u les constantes implicites ne d\'ependent que de $m$ et du
degr\'e $\deg(\delta)$ de $\delta$ dans un
plongement projectif.

  Par cons\'equent on obtient
\begin{equation}\label{eq2d}
\#\delta(\bbf_{\frp})=sq_{\frp}^m+O(q_{\frp}^{m-1/2}).
\end{equation}
\end{lemma}

\begin{proof}[\Dem]
Le premier cas de (\ref{eq2c}) d\'ecoule d'une version affaiblie du
Th\'eor\`eme de Lang-Weil \cite[Theorem
1]{lawe}. Le deuxi\`eme cas provient  du fait que les points
$\bbf_{\frp}$-rationnels de
$\delta_i$, pour $i>s$, sont contenus dans l'intersection des
composantes irr\'eductibles
conjugu\'ees de $\delta_i$ par
$G_{\bbf_{\frp}}$. La formule (\ref{eq2d}) suit imm\'edia\-tement de
(\ref{eq2c}).
\end{proof}

Nous allons  utiliser ce lemme pour estimer la contribution des
fibres singuli\`eres de notre
fibration. Notons $\Delta$
  le lieu discriminant de la fibration $f:\cax\to\cay$, c'est-\`a-dire~:
$$\Delta:=\left\{y\in\cay\;|\;\cax_y\hbox{ est singulier}\right\}.$$
Pour presque tout id\'eal premier $\frp$, la fibration
$\tilde{f}_{\frp}:\tilde{\cax}_{\frp}\rightarrow\tilde{\cay}_{\frp}$
est une fibration de vari\'et\'es lisses et a
pour lieu discriminant la r\'eduction  $\tilde{\Delta}_{\frp}$ de
$\Delta$ modulo $\frp$. Notons
$\tilde{\caf}_{\frp}=\ns_{\text{ver}}(\tilde{\cax}_{\frp})/\tilde{f}^*
_{\frp}(\ns(\tilde{\cay}_{\frp}))$, alors on
a, pour presque tout $\frp$, l'\'egalit\'e
$\tr(F_{\frp}|\tilde{\caf}_{\frp})=\tr(\Frob_{\frp}|{\caf}^{I_{\frp}})
$, ce qui permet le calcul
suivant.

\begin{proposition}\label{prop2a}
Pour presque tout id\'eal premier $\frp$ de $\cao_k$, on a~:
\begin{equation}\label{sfc}
\sum_{y\in\tilde{\Delta}(\bbf_{\frp})}\#\tilde{\cax}_{\frp,y}(\bbf_{\frp})
=q_{\frp}^{n-1}\tr\left(\Frob_{\frp}|\caf^{I_{\frp}}\right)+
q_{\frp}^{n-m}\#\tilde{\Delta}_{\frp}(\bbf_{\frp})+O(q_{\frp}^{n-3/2}).
\end{equation}
\end{proposition}

\begin{proof}[\Dem]
Posons $\cad:=f^{-1}(\Delta)=\cad_1\cup\dots\cad_{u}$ (o\`u les
$\cad_i$ d\'esignent les composantes
irr\'eductibles). Pour presque tout $\frp$ on a
$\tilde{\cad}_{\frp}=\tilde{f}^{-1}_{\frp}(\tilde{\Delta}_{\frp})$.
Le membre de gauche de la
formule (\ref{sfc}) est \'egal \`a $\#\tilde{\cad}_{\frp}(\bbf_{\frp})$.

Observons que
$\tr(\Frob_{\frp}|\caf^{I_{\frp}})=\tr(F_{\frp}|\tilde{\caf}_{\frp})$
est le nombre de g\'en\'erateurs
de
$\tilde{\caf}_{\frp}$ d\'efinis sur $\bbf_{\frp}$. Si $s$ (resp. $t$)
est le nombre de composantes
de
$\tilde{\cad}_{\frp}$ (resp. de $\tilde{\Delta}_{\frp}$) d\'efinies
sur $\bbf_{\frp}$, on a, d'apr\`es (\ref{eq2d})
du Lemme \ref{lem2b}~:
\begin{equation}\label{rd1}
\#\tilde{\cad}_{\frp}(\bbf_{\frp})=sq_{\frp}^{n-1}+O(q_{\frp}^{n-3/2})
\text{ et }
\#\tilde{\Delta}_{\frp}(\bbf_{\frp})=tq_{\frp}^{m-1}+O(q_{\frp}^{m-3/2}).
\end{equation}
Il reste \`a observer que $s=t+\tr(\Frob_{\frp}|\caf^{I_{\frp}})$
pour en d\'eduire la proposition (cf.
(\ref{eqns})).
\end{proof}

\begin{remark} Dans les travaux ant\'erieurs sur les fibrations en courbes (cf.
\cite{rosi},
\cite{wa} et
\cite{hp}), on d\'eterminait la contribution du nombre des points
rationnels dans
chaque fibre singuli\`ere individuellement et  on additionnait ensuite ces
contributions. Dans les deux premiers articles cit\'es, les fibres sont des
courbes de genre arithm\'etique un, et l'estimation a \'et\'e faite au cas par
cas en utilisant l'algorithme de Tate (d\'etaill\'e par exemple dans
\cite{siAEC}). Dans le troisi\`eme article cit\'e, les fibres sont des courbes
de genre arithm\'etique
$g$ et (d'apr\`es une remarque du referee - dans une premi\`ere
approche on utilisait
le travail de Artin-Winters {\cite{arwi}} pour obtenir une expression avec un
terme d'erreur de $O(q_{\frp}^{1/2})$) - on a d\'etermin\'e le nombre de points
rationnels des fibres singuli\`eres au moyen de la cohomologie $\ell$-adique
\`a support propre (cf. {\cite[Lemme 3.2]{hp}}).  Dans le cas
pr\'esent, au lieu
d'analyser chaque fibre, la proposition donne un contr\^ole de la somme totale
des contributions du nombre des points rationnels des fibres singuli\`eres.
Mais en employant la m\'ethode de d\'ecompte, on peut d\'emontrer la  g\'en\'eralisation suivante de
{\cite[Theorem
5.1]{wa}}.
\end{remark}

\begin{lemma}\label{lempts}
Pour chaque $y\in\tilde{\Delta}_{\frp}(\bbf_{\frp})$
soit $m_y$ le nombre des composantes $\bbf_{\frp}$-rationnelles de
$\tilde{f}_{\frp}^{-1}(y)$, alors on a
$$
\sum_{y\in\tilde{\Delta}_{\frp}(\bbf_{\frp})}(m_y-1)=q_{\frp}^{m-1}\tr(\Frob
_{\frp}\,|\,\caf^{I_{\frp}})+O(q_{\frp}^{m-3/2}).
$$ 
\end{lemma}

\begin{proof} 
Il suffit de combiner la formule de la proposition (\ref{prop2a}) avec l'\'egalit\'e suivante fournie par le lemme 
(\ref{lem2b})
$$
\#\tilde{\cax}_{\frp,y}(\bbf_{\frp})=m_yq_{\frp}^{n-m}+O(q_{\frp}^{n-m-1/2})
$$ 
En effet on en tire que
$$
\sum_{y\in\tilde{\Delta}_{\frp}(\bbf_{\frp})}\#\tilde{\cax}_{\frp,y}(\bbf_{\frp})=q_{\frp}^{n-m}\#\tilde{\Delta}_{\frp}
(\bbf_{\frp})+q_{\frp}^{n-m}\sum_{y\in\tilde{\Delta}_{\frp}(\bbf_{\frp})}(m_y-1)+O(q_{\frp}^{n-{3/2}})
$$
et le lemme suit.
\end{proof}




\section{Calculs de r\'esidus}

Dans ce paragraphe nous allons d\'emontrer le th\'eor\`eme suivant, dont le
Th\'eor\`eme \ref{th1a} est une cons\'equence imm\'ediate.

\begin{theorem}\label{th4ia}
Soit $f:\cax\to\cay$ une fibration de vari\'et\'es projectives
lisses, d\'efinies sur un corps de
nombres $k$, alors la fonction
\begin{multline}\label{eq4ia}
\sum_{\frp\notin
R}-\frA_{\frp}^*(\cax)\frac{\log(q_{\frp})}{q_{\frp}^s}+
\sum_{\frp\notin
R}\frB_{\frp}(\cax)\frac{\log(q_{\frp})}{q_{\frp}^{s+1}}\\
- \sum_{\frp\notin
R}b_{\frp}(\cax)\frac{\log(q_{\frp})}{q_{\frp}^{s+1}}+
\sum_{\frp\notin
R}b_{\frp}(\cay)\frac{\log(q_{\frp})}{q_{\frp}^{s+1}},
\end{multline}
qui est a priori analytique sur $\Re(s)>1$, se prolonge
analytiquement au demi-plan ferm\'e $\Re(s)\geq 1$ hormis
un p\^ole simple en $s=1$, avec un r\'esidu  \'egal \`a~:
\begin{equation}\label{eq4ib}
\rg\left(\frac{P_X(K)}{\tau
B(k)}\right)+\rg\left(\ns(X/K)\right)-\rg\left(\ns(\cax/k)\right)
+\rg\left(\ns(\cay/k)\right).
\end{equation}
\end{theorem}

Commen\c{c}ons par rappeler que si
$$
L(s):=\prod_{\frp}\det(1-\Frob_{\frp}q_{\frp}^{-s}|V^{I_{\frp}})^{-1}
$$
est une fonction $L$ de poids $w$ (i.e., les valeurs propres de
$\Frob_{\frp}$ sont de modules $q_{\frp}^{w/2}$),
alors
$$
-\frac{L'(s)}{L(s)}=\frac d{ds}(\log
L(s))=\sum_{\frp}\tr(\Frob_{\frp}|V^{I_{\frp}})\frac{\log(q_{\frp})}{q
_{\frp}^s}+g(s)
$$
avec $g(s)$ holomorphe sur le demi-plan $\Re(s)>(w+1)/2$. Dans le cas
o\`u le groupe de Galois
agit \`a travers un groupe fini, on dispose du th\'eor\`eme classique suivant.

\begin{proposition}[Artin, Brauer]\label{prop3b}\cite[Proposition 1.5.1]{rosi}
Soit $V$ un $\bbq$-espace vecto\-riel de dimension finie avec une
action continue de $G_k$, et soit $L(V,s)$
la fonction $L$ d'Artin associ\'ee.

\begin{enumerate}
\item $L(V,s)$ a une continuation m\'eromorphe \`a $\bbc$.

\item $L(V,s)$ est holomorphe sur la droite $\Re(s)=1$ \`a
l'exception \'eventuelle du point $s=1$, o\`u l'on
a
$\ord_{s=1}(L(V,s))=-\dim(V^{G_k})$.

\item $L(V,s)$ ne s'annule pas sur la droite $\Re(s)=1$.
\end{enumerate}
\end{proposition}

En particulier, en choisissant $V=\caf\otimes\bbq$,
on obtient le comportement suivant~:
\begin{equation}\label{4a}
\res_{s=1}\left(\sum_{\frp}\tr(\Frob_{\frp}|\caf^{I_{\frp}})\frac{\log
{q_{\frp}}} {q_{\frp}^s}
\right)=\rg(\caf^{G_k})
\end{equation}

\begin{proof}[\Dem\ du Th\'eor\`eme \ref{th4ia}]
La formule de Lefschetz (cf. Lemme \ref{lef}) appli\-qu\'ee \`a $\cax$
et $\cay$ permet d'\'ecrire
que,
\begin{equation}\label{uu00}
\#\tilde{\cax}_{\frp}(\bbf_{\frp})=q_{\frp}^{n}-a_{\frp}(\cax)q_{\frp}^{n-1}+
b_{\frp}(\cax)q_{\frp}^{n-2}+O(q_{\frp}^{n-3/2}),
\end{equation}

\begin{equation}\label{uu01}
\#\tilde{\cay}_{\frp}(\bbf_{\frp}) =q_{\frp}^{m}-a_{\frp}(\cay)q_{\frp}^{m-1}+
b_{\frp}(\cay)q_{\frp}^{m-2}+O(q_{\frp}^{m-3/2}).
\end{equation}

La formule de Lefschetz (cf. Lemme \ref{lef}) appliqu\'ee \`a
$\tilde{\cax}_{\frp,y}$ permet d'\'ecrire que,
pour
$y\in(\tilde{\cay}_{\frp}-\tilde{\Delta}_{\frp})(\bbf_{\frp})$, on a
\begin{equation}\label{uu2}
\#\tilde{\cax}_{\frp,y}(\bbf_{\frp})=q_{\frp}^{n-m}-a_{\frp}(\tilde{\cax}_{\frp,y})q_{\frp}^{n-m-1}+
b_{\frp}(\tilde{\cax}_{\frp,y})q_{\frp}^{n-m-2}+O(q_{\frp}^{n-m-3/2}).
\end{equation}
On peut donc \'ecrire en invoquant   (\ref{uu2}) et la Proposition \ref{sfc}

\begin{equation}\begin{aligned}\label{uu3}
&\#\tilde{\cax}_{\frp}(\bbf_{\frp})=
\sum_{y\in(\tilde{\cay}_{\frp}-\tilde{\Delta}_{\frp})(\bbf_{\frp})}
\#\tilde{\cax}_{\frp,y}(\bbf_{\frp})+\sum_{y\in\tilde{\Delta}_{\frp}(\
bbf_{\frp})}\#\tilde{\cax}_{\frp,y}(\bbf_{\frp})\\
&=\sum_{y\in(\tilde{\cay}_{\frp}-\tilde{\Delta}_{\frp})(\bbf_{\frp})}
(q_{\frp}^{n-m}-a_{\frp}(\tilde{\cax}_{\frp,y})q_{\frp}^{n-m-1}+
b_{\frp}(\tilde{\cax}_{\frp,y})q_{\frp}^{n-m-2}\\
&+O(q_{\frp}^{n-m-3/2}))+q_{\frp}^{n-1}\tr(\Frob_{\frp}|\caf^{I_{\frp}})+
q_{\frp}^{n-m}\#\tilde{\Delta}_{\frp}(\bbf_{\frp})+O(q_{\frp}^{n-3/2})\\
&=q_{\frp}^{n-m}\#\tilde{\cay}_{\frp}(\bbf_{\frp})-q_{\frp}^{n-m-1}
\left(\sum_{y\in\tilde{\cay}_{\frp}(\bbf_{\frp})}a_{\frp}(\tilde{\cax}
_{\frp,y})-
\sum_{y\in\tilde{\Delta}_{\frp}(\bbf_{\frp})}a_{\frp}(\tilde{\cax}_{\frp,y})\right)\\
&+q_{\frp}^{n-m-2}\left(\sum_{y\in\tilde{\cay}_{\frp}(\bbf_{\frp})}b_{ \frp}(\tilde{\cax}_{\frp,y})-
\sum_{y\in\tilde{\Delta}_{\frp}(\bbf_{\frp})}b_{\frp}(\tilde{\cax}_{\frp,y})\right)+
q_{\frp}^{n-1}\tr(\Frob_{\frp}\,|\,\caf^{I_{\frp}})\\
&+O(q_{\frp}^{n-3/2})=q_{\frp}^{n}-a_{\frp}(\cay)q_{\frp}^{n-1}-\frA_{
\frp}(\cax)q_{\frp}^{n-1}
+b_{\frp}(\cay)q_{\frp}^{n-2}\\
& +\frB_{\frp}(\cax)q_{\frp}^{n-2}+q_{\frp}^{n-1}\tr(\Frob_{\frp}\,|\,\caf^{I_{\frp}})+O(q_{\frp}^{n-3/2}).
\end{aligned}\end{equation}

Pour la derni\`ere \'egalit\'e, on a utilis\'e la d\'efinition des
$\frA_{\frp}(\cax)$ et $\frB_{\frp}(\cax)$ ainsi
que les estimations
\begin{equation}\label{eqweil2}
\sum_{y\in\tilde{\Delta}_{\frp}(\bbf_{\frp})}a_{\frp}(\tilde{\cax}_{\frp,y})=O(q_{\frp}^{m-1/2})\text{ et }
\sum_{y\in\tilde{\Delta}_{\frp}(\bbf_{\frp})}b_{\frp}(\tilde{\cax}_{\frp,y})=O(q_{\frp}^{m})
\end{equation}
pour faire rentrer ces deux termes dans le terme reste. Pour obtenir
les estimations de (\ref{eqweil2}) on a
employ\'e un r\'esultat d\^u \`a Deligne \cite[Theorem 3.3.1]{del81} qui dit
que les valeurs propres de $F_{\frp}$ agissant sur
$H^1_c(\tilde{\cax}_{\frp,y})$ (resp. $H^2_c(\tilde{\cax}_{\frp,y})$) sont de
valeur absolue au plus
$q_{\frp}^{1/2}$ (resp. $q_{\frp}$). Finalement on emploie les
estimations de Lang et Weil (cf.
Lemme \ref{lem2b}) pour conclure (\ref{eqweil2}).

En invoquant la d\'efinition des  $\frA^*_{\frp}(\cax)$, l'\'egalit\'e
$a_{\frp}(\cax)=a_{\frp}(\cay)+a_{\frp}(B)$ (voir la Proposition
\ref{prop2c}) et (\ref{uu00}), on
d\'eduit de l'\'egalit\'e (\ref{uu3}) (divis\'ee par
$q_{\frp}^{n-1}$) que
\begin{equation}\label{keyf}
-\frA^*_{\frp}(\cax)+\frac{\frB_{\frp}(\cax)}{q_{\frp}}-\frac{b_{\frp}(\cax)}{
q_{\frp}}+\frac{b_{\frp}(\cay)}{q_{\frp}}=-\tr(\Frob_{\frp}\,|\,\caf^{
I_{\frp}})+O(q_{\frp}^{-1/2}).
\end{equation}
D'apr\`es la Proposition \ref{prop3b} et en particulier (\ref{4a}), on obtient
\begin{equation}\begin{aligned}\label{fin}
&\res_{s=1}\left(\sum_{\frp\notin R}\left(-\frA_{\frp}^*(\cax)+
\frac{\frB_{\frp}(\cax)}{q_{\frp}}-
\frac{b_{\frp}(\cax)}{q_{\frp}}+ \frac{b_{\frp}(\cay)}{q_{\frp}}
\right)\frac{\log(q_{\frp})}{q_{\frp}^s}\right)\\
&=-\rg\left(\caf^{G_k}\right)= \rg\left(\frac{P_X(K)}{\tau
B(k)}\right)+\rg\left(\ns(X/K)\right)\\
&-\rg\left(\ns(\cax/k)\right) +\rg\left(\ns(\cay/k)\right),
\end{aligned}\end{equation}
o\`u la derni\`ere \'egalit\'e suit de la formule de Shioda-Tate
(Proposition \ref{prop2b}). Ce qui ach\`eve la
preuve du Th\'eor\`eme \ref{th4ia}.
\end{proof}




\section{Illustrations}

On peut facilement obtenir des exemples de fibrations gr\^ace \`a la
proposition classique
suivante.

\begin{proposition}\label{prop5a}
Soit $\cax$ une vari\'et\'e lisse et projective de dimension $n$,
d\'efinie sur un corps infini $k$. Quitte
\`a remplacer $\cax$ par son \'eclatement $\cax'$ le long d'une
sous-vari\'et\'e lisse de
dimension $n-m-1$, il existe une fibration en vari\'et\'es de
dimension $n-m$, disons
$f:\cax'\to\bbp^{m}$, \'egalement d\'efinie sur
$k$.
\end{proposition}

\begin{proof}[\Dem] (Cf \cite[Proposition 5.1]{hp} pour le cas d'une
fibration en courbes).
Choisissons $\phi:\cax\rightarrow\bbp^n$ un morphisme fini (par
exemple obtenu par normalisation de Noether, i.e.,
par un plongement dans un $\bbp^N$ suivi d'une projection lin\'eaire)
et choisissons une
sous-vari\'et\'e lin\'eaire g\'en\'erale $L$ d\'efinie sur $k$, de
dimension $n-m-1$ dans
$\bbp^n$. Alors
$\phi^{-1}(L)$ est lisse (et m\^eme connexe si
$n-m>1$). Notons ${\bbp}'$ l'\'eclat\'e de $\bbp^n$ au-dessus de $L$
et par $\pi':\bbp'\to\bbp^m$
la fibration en espaces projectifs de dimension $n-m$ qui s'en d\'eduit. Soit
$\pi:\cax'\to\cax$ l'\'eclatement de $\cax$ au-dessus de
$\phi^{-1}(L)$, on a donc un morphisme fini $\phi':\cax'\to{\bbp}'$
qui, compos\'e avec $\pi'$, fournit la
fibration cherch\'ee $f:=\phi'\circ\pi':\cax'\to\bbp^m$  en
vari\'et\'es de dimension $n-m$.
\end{proof}

Pour toute vari\'et\'e projective lisse d\'efinie sur un corps de
nombres $k$, d\'enotons
$H^2(\cax)(1):=H^2_{\et}(\cax\times_k\ov{k},\bbq_{\ell}(1))$. Soit
$\text{cl}_{\ell}:\ns(\cax)\to
H^2(\cax)(1)$ l'application de classes de diviseurs.

La proposition suivante montre que la Conjecture T est invariante par
une application $k$-birationnelle qui
est compos\'ee d'un nombre fini d'\'eclatements.

\begin{proposition}\label{prop5b}
Soit $\phi:\cax'\to\cax$ une application $k$-birationnelle entre deux
vari\'et\'es lisses et projectives de
dimension
$n$ d\'efinies sur un corps de nombres $k$. Soit
$\phi':\cax'\times_k\ov{k}\to\cax\times_k\ov{k}$
l'application
$\ov{k}$-birationnelle d\'eduite de $\phi$ par extension de
scalaires. Supposons que $\phi'$ soit compos\'ee
d'un nombre fini d'\'eclatements de centres lisses. Donc, la
Conjecture T est vraie pour $\cax'/k$
si et seulement si elle est vraie pour $\cax/k$.
\end{proposition}

\begin{proof}[\Dem]
Soit
$\cax\times_k\ov{k}=\cax_0\overset{\phi_1}\longleftarrow\cax_1\overset
{\phi_2}\longleftarrow\cdots
\overset{\phi_{r-1}}\longleftarrow\cax_{r-1}\overset{\phi_r}
\longleftarrow\cax_r=\cax'\times_k\ov{k}$ la factorisation de
$\phi'$, o\`u les applications
$\phi_i$ sont des \'eclatements de centres lisses. Chacune de ces
transformations rajoute un
diviseur exceptionnel $E_i$ au groupe de N\'eron-Severi et au groupe
de cohomologie,
c'est-\`a-dire, soit $\{E_i\}$ la classe de $E_i$ dans
$\ns(\cax_i)$, donc
\begin{multline*}
\ns(\cax_i)\cong\phi_i^*(\ns(\cax_{i-1}))\oplus\bbz[\{E_i\}]\text{ et}\\
H^2(\cax_i)(1)\cong\phi_i^*(H^2(\cax_{i-1})(1))\oplus\bbq_{\ell}(1)[\text{cl}_{\ell}(\{E_i\})].
\end{multline*}
Globalement on obtient donc, en notant par $V$ le groupe des classes
de diviseurs engendr\'e par les
diviseurs exceptionnels sur $\cax'\times_k\ov{k}$, et remarquant que
la restriction de
$\text{cl}_{\ell}$ \`a
$V\otimes\bbq_{\ell}$ donne un isomorphisme sur son image $\mathcal{V}$,
$$
\ns(\cax')\cong{\phi'}^*(\ns(\cax))\oplus V\text{ et }
H^2(\cax')(1)\cong{\phi'}^*(H^2(\cax)(1))\oplus (V\otimes\bbq_{\ell})
$$
Le quotient
$$
\frac{L_2(\cax'/k,s)}{L_2(\cax/k,s)}=\frac{L(H^2(\cax')(1),s-1)}{L(H^2
(\cax)(1),s-1)}=L(V,s-1)
$$
est donc une fonction $L$ d'Artin dont l'ordre en $s=2$ (cf.
Proposition \ref{prop3b}) est \'egal
\`a
\begin{equation}\label{eq5a}\begin{aligned}
-\ord_{s=2}(L(V,s-1))&=\dim((V\otimes\bbq_{\ell})^{G_k})\\
&=\dim(H^2(\cax')(1)^{G_k})- \dim(H^2(\cax)(1)^{G_k}).
\end{aligned}\end{equation}

D'autre part on a bien
$$\rg(\ns(\cax'/k))-\rg(\ns(\cax/k))
=\dim((V\otimes\bbq_{\ell})^{G_k})$$ d'o\`u le r\'esultat.
\end{proof}

\begin{remark}
Au contraire du cas des surfaces, on ne peut pas dire qu'une
application $\ov{k}$-birationnelle entre
vari\'et\'es projectives et lisses de dimension $n\ge3$ sur $\ov{k}$
soit compos\'ee par un nombre
fini d'\'eclatements. Cependant, le (profond) Th\'eor\`eme de
Factorisation Faible (cf.
\cite{wl} et \cite{ab}) affirme que, avec les notations ci-dessus, si
$\phi':\cax'\times_k\ov{k}\to\cax\times_k\ov{k}$ est une application
$\ov{k}$-birationnelle, il existe
une factorisation de $\phi'$,
$\cax\times_k\ov{k}=\cax_0\overset{\phi_1}\longleftarrow\cax_1\overset
{\phi_2}\longleftarrow\cdots
\overset{\phi_{r-1}}\longleftarrow\cax_{r-1}\overset{\phi_r}
\longleftarrow\cax_r=\cax'\times_k\ov{k}$, o\`u chaque
$\phi_i$ est un \'eclatement ou une contraction de centre lisse. Le
cas de la factorisation forte, d'o\`u
d\'ecoulerait la Proposition \ref{prop5b}, reste encore un probl\`eme ouvert.
\end{remark}

\begin{remark}\label{lestates}
La Conjecture T est en effet cons\'equence de deux conjectures (qui sont
prouv\'ees dans tous les cas o\`u la validit\'e de cette conjecture est connue,
cf. \cite{ram}, \cite{rosi}):

\begin{itemize}
\item (\textbf{Conjecture de Tate I}) \emph{On a un isomorphisme
de $\bbq_{\ell}$-espaces
vectoriels}
\begin{equation}\label{eq5b}
\ns(\cax/k)\otimes\bbq_{\ell}\cong H^2({\cax})(1)^{G_k}
\end{equation}

\item (\textbf{Conjecture de Tate II}) \emph{On a l'\'egalit\'e}
\begin{equation}\label{eq5c}
-\ord_{s=2}(L_2(\cax/k,s))=\dim(H^2({\cax})(1)^{G_k}).
\end{equation}

\end{itemize}

\end{remark}

\begin{remark}\label{rem5a}
Lorsque les classes de diviseurs engendrent $H^2(\cax)(1)$, alors la
Conjecture T est vraie pour $\cax$ (la
Conjecture I est triviale et la Conjecture II d\'ecoule de la
Proposition \ref{prop3b}). En
particulier la Conjecture T est vraie pour $\bbp^n$ et pour toute
vari\'et\'e birationnelle \`a
$\bbp^n$.
\end{remark}

Nous terminerons cette section en notant des cons\'equences du Th\'eor\`eme
\ref{th1a} analogues \`a celles donn\'ees dans
\cite{hp}. Dans le reste de cette section, aussi bien que dans la
prochaine, dire que les Conjectures T$_{\text{fin}}$, M$_{\text{an}}$ et
BSD$_{\text{fin}}$ sont vraies pour les fibrations en vari\'et\'es
veut dire que
ces conjectures sont vraies pour toute $d$-fibration pour tout $d\ge1$ entier.

\begin{theorem}\label{th5a} Les Conjectures suivantes sont \'equivalentes:

\begin{enumerate}
\item[(a)] La Conjecture T.

\item[(b)] La Conjecture M$_{\text{an}}$ pour les fibrations en vari\'et\'es.

\item[(c)] La Conjecture N$_{\text{an}}$ pour les fibrations en courbes.

\item[(d)] La Conjecture M$_{\text{an}}$ pour les fibrations en vari\'et\'es au
dessus d'un espace projectif.

\item[(e)] La Conjecture N$_{\text{an}}$ pour les fibrations en courbes au
dessus d'un espace projectif.

\end{enumerate}

\end{theorem}

\begin{proof}[\Dem]

Nous d\'emontrerons que \mbox{$ a) \Rightarrow b) \Rightarrow c)
\Rightarrow e) \Rightarrow a)$}.  Le th\'eor\`eme en d\'ecoule, puisqu'il est
facile \`a voir que \mbox{$ b) \Rightarrow d) \Rightarrow e) $}.

Soit $f:\cax\to \cay$ une $d$-fibration en vari\'et\'es telle que
$f$, $\cax$ et
$\cay$ soient d\'efinis sur $k$ avec
$\dim(\cax)=n$, $\dim(\cay)=m$ et $d=n-m>0$. Par le Th\'eor\`eme \ref{th1a}, la
validit\'e de la Conjecture T pour
$\cax/k$ et $\cay/k$ entra\^{\i}ne la Conjecture M$_{\text{an}}$ pour $f$.

Il est trivial que la Conjecture M$_{\text{an}}$ pour les fibrations en
  vari\'et\'es implique la Conjecture M$_{\text{an}}$ pour les fibrations en
courbes.  Comme nous l'avons d\'ej\`a observ\'e dans Remarque {\ref{remnm}}, la
Conjecture M$_{\text{an}}$ pour les fibrations en courbes est \'equivalente \`a
la Conjecture N$_{\text{an}}$ pour les fibrations en courbes.  En outre, il est
clair que la Conjecture N$_{\text{an}}$ pour les fibrations en
courbes au dessus
d'une vari\'et\'e quelconque implique la Conjecture N$_{\text{an}}$ pour les
fibrations en courbes au dessus d'un espace projectif.

Supposons maintenant vraie la Conjecture N$_{\text{an}}$ pour les
fibrations en courbes au dessus d'un espace projectif. Soit $\cax/k$
une vari\'et\'e lisse projective de
dimension
$n$. Par la Proposition \ref{prop5a} (o\`u bien \cite[Proposition
5.1]{hp}), $\cax$ est
$k$-birationnellement \'equivalente \`a une vari\'et\'e
${\cax'}$ qui se fibre en courbes sur $\bbp^{n-1}$. Par la Remarque
\ref{rem5a}, la conjecture de
Tate est vraie pour
$\bbp^{n-1}$, et il suit du Th\'eor\`eme \ref{th1a} que la conjecture
de Tate est vraie
pour ${\cax'}$, donc, par la
Proposition \ref{prop5b}, elle est aussi vraie pour $\cax$.
\end{proof}

\begin{remark}\label{rem5b}
Dans la d\'emonstration du Th\'eor\`eme \ref{th5a}, il suffit
d'utiliser \cite[Proposition 5.1]{hp} pour obtenir une fibration en
courbes, parce que la preuve que la Conjecture M$_{\text{an}}$
  pour les fibrations en vari\'et\'es implique la Conjecture T est indirecte
(\`a travers la Conecture N$_{\text{an}}$ pour les fibrations en courbes); mais
en utilisant la Proposition \ref{prop5a} on pourrait d\'emontrer directement
l'\'equivalence entre la Conjecture T et la Conjecture M$_{\text{an}}$.
\end{remark}




\section{Les rapports entre les Conjectures T$_{\text{fin}}$,
N$_{\text{an}}$ et BSD$_{\text{fin}}$
pour les fibrations en vari\'et\'es}

\begin{proposition}\label{prop7a}
La Conjecture T$_{\text{fin}}$ pour les fibrations en vari\'et\'es
entra\^{\i}ne la Conjecture T.
\end{proposition}

\begin{proof}[\Dem]

Soit $\cax/k$ une vari\'et\'e projective lisse de dimension $n$ d\'efinie sur
corps de nombres $k$, $m\ge1$ un nombre entier, et $\cax'=\cax\times\bbp^m$.
Soit $p_1:\cax'\to\cax$ la projection sur le premier facteur, et
$f:\cax'\to\bbp^m$ la $n$-fibration \'egale \`a la projection sur le deuxi\`eme
facteur.  Allors $f$ est une fibration constante, o\`u toutes les fibres sont
\'egales \`a $\cax$. D\'enotons $\bbp^m$ par
$\cay$. Donc, pour tout $\frp\notin R$ et
$y\in\tilde{\cay}_{\frp}(\bbf_{\frp})$ on a

\begin{equation*}
b_{\frp}(\tilde{\cax}'_{\frp,y})=b_{\frp}(\tilde{\cax}_{\frp})=\tr(F_{
\frp}\,|\,
H^2(\tilde{\cax}_{\frp}))
=\tr(\Frob_{\frp}\,|\,H^2(\cax)^{I_{\frp}})=b_{\frp}(\cax).
\end{equation*}
D'o\`u,

$$
\frB_{\frp}(\cax')=\frac1{q_{\frp}^m}\sum_{y\in\tilde{\cay}_{\frp}(\bbf_{\frp})} b_{\frp}(\tilde{\cax}_{\frp,y}')
=\left(1+\frac1{q_{\frp}}+\ldots+\frac1{q_{\frp}^m}\right)b_{\frp}(\cax).
$$

Le calcul des r\'esidus nous donne
\begin{equation}\label{eq7a}
\res_{s=2}\left(\sum_{\frp\notin
R}\frB_{\frp}(\cax')\frac{\log(q_{\frp})}{q_{\frp}^s}\right)=
\res_{s=2}\left(\sum_{\frp\notin
R}b_{\frp}(\cax)\frac{\log(q_{\frp})}{q_{\frp}^s}\right).
\end{equation}

Soit $k(\un{t}):=k(t_1,\cdots,t_m)$ le corps de fonctions de
$\bbp^m_k$ et $\cax'_{\eta}/k(\un{t})$ la
fibre g\'en\'erique de $f$. Le r\'esultat est donc \'equivalent \`a prouver
\begin{equation}\label{eq7b}
\rg(\ns(\cax'_{\eta}/k(\un{t})))=\rg(\ns(\cax/k)).
\end{equation}
Mais $\cax'_{\eta}\cong\cax\times_kk(\un{t})$, donc
$\ns(\cax'_{\eta}/k(\un{t}))\cong\ns(\cax/k)$, d'o\`u
l'\'egalit\'e entre les rangs suit.
\end{proof}

\begin{theorem}\label{th7b}
Les Conjectures N$_{\text{an}}$ et T$_{\text{fin}}$ pour les
fibrations en vari\'et\'es sont
\'equivalentes.
\end{theorem}

\begin{proof}[\Dem]
Supposons que la Conjecture T$_{\text{fin}}$ soit vraie pour les
fibrations en vari\'et\'es. Soit $f:\cax\to\cay$
une fibration en vari\'et\'es d\'efinie sur un corps de nombres $k$.
Par la Proposition
\ref{prop7a}, la Conjecture T est vraie pour $\cax/k$ et pour
$\cay/k$. Ainsi, par le Corollaire
\ref{co1a}, la Conjecture N$_{\text{an}}$ est vraie pour
$f$.

R\'eciproquement, supposons que la Conjecture N$_{\text{an}}$ soit
vraie pour les fibrations en vari\'et\'es,
a fortiori elle est vraie pour les fibrations en courbes. Soit
$f:\cax\to\cay$ une fibration en
vari\'et\'es. Par le Th\'eor\`eme \ref{th5a}, la Conjecture T est vraie pour
$\cax/k$ et $\cay/k$. Encore une fois, on conclut par le Corollaire
\ref{co1a} que la Conjecture T$_{\text{fin}}$ est vraie pour $f$.
\end{proof}

\begin{proposition}\label{prop7b}
La Conjecture BSD$_{\text{fin}}$ pour les fibrations  en vari\'et\'es
entra\^{\i}ne la Conjecture de Birch
\& Swinnerton-Dyer pour les vari\'et\'es ab\'eliennes sur un corps de nombres.
\end{proposition}

\begin{proof}[\Dem]

Soit $k$ un corps de nombres, $A/k$ une vari\'et\'e ab\'elienne de
dimension $d$, $m\ge1$ un nombre
entier,
$\caa=A\times\bbp^m$ et $f:\caa\to\bbp^m$ la $d$-fibration constante
obtenue par la projection dans la
deuxi\`eme composante. Soit $k(\un{t}):=k(t_1,\cdots,t_m)$ le corps
de fonctions de $\bbp^m_k$.
Par construction, la
$k(\un{t})/k$-trace de la fibre g\'en\'erique $\caa_{\eta}/k(\un{t})$
de $f$ est \'egale \`a $A$. De plus,
comme
$\caa_{\eta}\cong A\times_kk(\un{t})$, on a
$\caa_{\eta}(k(\un{t}))\cong A(k(\un{t}))=A(k)$. Alors,
\begin{equation}\label{eq7ic}
\rg(\caa_{\eta}(k(t)))=\rg(A(k)).
\end{equation}

D\'enotons $\bbp^m$ par $\cay$. Pour tout $\frp\notin R$ et
$y\in\tilde{\cay}_{\frp}(\bbf_{\frp})$ on
a
\begin{equation*}
a_{\frp}(\tilde{\caa}_{\frp,y})=a_{\frp}(\tilde{A}_{\frp})=\tr(F_{\frp
}\,|\,H^1(\tilde{A}_{\frp}))
=\tr(\Frob_{\frp}\,|\,H^1(A)^{I_{\frp}})=a_{\frp}(A).
\end{equation*}
D'o\`u,

$$
\frA_{\frp}(\caa)=\frac1{q_{\frp}^m}\sum_{y\in\tilde{\cay}_{\frp}(\bbf_{\frp})}
a_{\frp}(\tilde{\caa}_{\frp,y})=
\left(1+\frac1{q_{\frp}}+\ldots+\frac1{q_{\frp}^m}\right)a_{\frp}(A).
$$
Le calcul de r\'esidus donne

$$
\res_{s=1}\left(\sum_{\frp\notin
R}-\frA_{\frp}(\caa)\frac{\log(q_{\frp})}{q_{\frp}^s}\right)=
\res_{s=1}\left(\sum_{\frp\notin
R}-a_{\frp}(A)\frac{\log(q_{\frp})}{q_{\frp}^s}\right).
$$
Le r\'esultat se d\'eduit de cette \'egalit\'e, de (\ref{eq7ic}) et
de (\ref{eqbsdfin}).
\end{proof}

\begin{proposition}\label{prop7c}
La Conjecture BSD$_{\text{fin}}$ pour les fibrations en vari\'et\'es
entra\^{\i}ne la Conjecture N$_{\text{an}}$
pour les fibrations en vari\'et\'es.
\end{proposition}

\begin{proof}[\Dem] Soit $f:\cax\to\cay$ une fibration en vari\'et\'es,
$P_{\cax_{\eta}}$ la vari\'et\'e de Picard de la fibre g\'en\'erique
$\cax_{\eta}/k(\cay)$ de $f$ et $(\tau,B)$ sa $k(\cay)/k$-trace. Par
hypoth\`ese, la Conjecture BSD$_{\text{fin}}$ est vraie pour les fibrations en
vari\'et\'es, donc par la Proposition \ref{prop7b}, la Conjecture de Birch et
Swinnerton-Dyer est vraie pour $B/k$. Notons enfin que la validit\'e de la
conjecture N$_{\text{an}}$ pour
$f$ est une cons\'equence de la Conjecture de Birch et Swinnerton-Dyer pour
$B/k$ et de la Conjecture BSD$_{\text{fin}}$ pour $f$.
\end{proof}

\begin{theorem}\label{th7c}
La Conjecture BSD$_{\text{fin}}$ pour les fibrations en vari\'et\'es
entra\^{\i}ne la Conjecture T$_{\text{fin}}$
pour les fibrations en vari\'et\'es.
\end{theorem}

\begin{proof}[\Dem]
Le th\'eor\`eme suit imm\'ediatement du Th\'eor\`eme \ref{th7b} et de
la Proposition \ref{prop7c}.
\end{proof}




\end{document}